\documentclass[12pt]{article}
\usepackage{amsmath,amsthm,amsfonts,amssymb,amscd}
\usepackage[usenames]{color}
\usepackage[dvips]{graphicx}
\usepackage{epsfig}
\def\phi{\varphi}

\begin{document}

\centerline{\bf ON THE PROBABILITY OF CO-PRIMALITY OF TWO NATURAL}
\centerline{\bf NUMBERS CHOSEN AT RANDOM}
\medskip
\centerline{\bf (Who was the first to pose and solve this problem?)}

\bigskip

\centerline{\it Sergei Abramovich$^{a}$, Yakov Yu. Nikitin$^{b, c}$}

\medskip

\centerline{$^{a}$ School of Education and Professional Studies, State University of New York, }
\centerline{Potsdam, NY 13676, USA, abramovs@potsdam.edu }

\medskip

\centerline { $^{b}$ Department of Mathematics and Mechanics, Saint Petersburg State University,}
\centerline{ 7/9 Universitetskaya nab., St. Petersburg, 199034, Russia, y.nikitin@spbu.ru}

\medskip

\centerline{ $^{c}$ National Research University - Higher School of Economics,}
\centerline{ Souza Pechatnikov, 16, St.Peters\-burg 190008, Russia}

\bigskip

  {\small The article attempts to demonstrate the rich history of one truly remarkable problem situated at the confluence of probability theory and theory of numbers - finding the probability of co-primality of two randomly selected natural numbers. The goal of the article is to reveal the genesis of the problem and understand the role of mathematicians of the past in solving this problem and using it as a background for advancing mathematical ideas. The article describes different approaches to the solution of the problem, reviews its various generalizations, examines similar problems, and offers diverse historical perspectives.}

\medskip


\bigskip

\section{ Introduction}

\indent	The problem of finding probability that two randomly chosen natural numbers are relatively prime (referred to below as the Problem) is one of the gems of mathematical heritage. Its equivalent formulation is finding the probability  of the irreducibility of a fraction with randomly selected numerator and denominator from the set of natural numbers. The Problem  is quite remarkable despite, or perhaps because of the ease of its formulation, when fractions, which are used as tools in computing chances of simple events already at the primary school level, are themselves considered through the lenses of probability.

Yet, the Problem's appeal goes far beyond the elementary level and it may safely be regarded as one of the most famous problems situated at the confluence of the theory of probability and analytic number theory. It can be found in many text/problem-solving books, articles, and monographs. Here are just two examples: in a famous textbook on number theory, Hardy and Wright (1975), the Problem is presented as Theorem 332; in a more recent book on analytic and probabilistic number theory it is a corollary from Theorem 4, see Tenenbaum (1995), p. 40; in mathematically beautiful and remarkable by originality problem-solving book by  Chaumont and Yor (2003) it is exercise 2.10.

	All this motivated the authors to address the following queries. What is the genesis of the Problem? How could it be traced back to the work of others? Who are those others? How can the problem be used to advance our knowledge of the history of mathematics? In this paper, different approaches to the Problem will be described, various generalizations reviewed, reflections on similar problems offered, and some historical perspectives inferred.

\section{ Three aproaches to solving the problem}

 There exist three main approaches to solving the Problem. The first two are based on analytic and algebraic number theories. The third one is based on probabilistic reasoning which progressed from being heuristic to becoming rigorous.

	From the number theory perspective, one can take a large natural number $n,$ then consider all possible pairs of numbers in the range 1 through $n$ $$
B_n = \{ {i,k} : 1 \le i < k \le n\}
$$
and select from them the set of relatively prime pairs, that is, the set
$$
 A_n = \{ {i,k} : 1 \le i < k \le n, GCD(i,k) =1\}.
$$
Consequently, the probability $P$ that two numbers selected at random are relatively prime is quite natural to define as the limit (assuming that it exists)
$$
P = \lim_{n \to \infty} P_n:= \lim_{n \to \infty}\frac{|A_n|}{|B_n|} = \lim_{n \to \infty}\frac{2|A_n|}{n(n-1)}.
$$

That is how the Problem can be solved through a number theory approach. It should be noted though, that the probability is substituted here with the {\it asymptotic density} of the set of irreducible pairs that, unlike probability, is not countably additive. Originally, in order to compute the asymptotic behavior of $|A_n|$ , some non-rigorous reasoning was used.

	Then, it became clear that $|A_n|$  is equal to the sum of values of Euler phi function $\sum_{k=2}^n \varphi(k)$  the asymptotic behavior of which was found first by Dirichlet (1849, 1897) and then (with a better residual term) by Mertens (1874), Kronecker's student, by using the M\"{o}bius function, namely,
$$
|A_n| \sim \frac{3}{\pi^2} n^2, \ n \to \infty,
$$
from where, finally, the well-known value $P = 6/\pi^2 \approx 0.607927...  \ $  results.

	Later, this solution, with gradual simplifications, has been repeated multiple times, e.g., in Sylvester (1883); Kronecker (1901); Hardy and Wright (1975); Apostol (1976); Knuth (1981).

	Another method of solving the Problem is based on heuristic reasoning when the key step is to use, yet without sufficient justification, the probabilistic property of "independence" of the events comprised of divisibility of  randomly selected natural number by different prime numbers. After a series of reasoning techniques and algebraic manipulations, the sought probability can be represented through the infinite product $\prod_{p} \left(1- \frac{1}{p^2}\right )$  over all prime numbers (see details in Yaglom and Yaglom (1964/67/87), or in Schroeder (2008), pp. 52-53). Due to the famous identity
$$
\prod_{p} \left(1- \frac{1}{p^2}\right )^{-1} = \sum_{n\ge1} \frac{1}{n^2} = \pi^2/6, \qquad  \qquad  \qquad    \qquad  \qquad           (1)	
$$					
discovered by Euler in 1737, once again, the same value $P$ as before results. That was the way many great mathematicians of the past had have approached and solved the Problem, although indicating that their reasoning was not quite rigorous. Even nowadays this solution is presented as a new one, see Hombas (2013).
	
Recall that the function
$$
\zeta(s) = \sum_{n\ge1} \frac{1}{n^s}
$$
defined for real and complex values of $s,$ is called zeta function. For real values of s, it was introduced in 1737 by Euler who was able to compute all the values $\zeta(2m)$  of the function for even $s.$ Notwithstanding, up to this date, the exact values of the function for odd $s$ are not known, although their approximate values may be easily computed to any given accuracy.

	In 1859, Riemann studied the zeta function  for complex values of $s$ and revealed the close relationship between the non-trivial zeroes of the function and the distribution of prime numbers. The Riemann hypothesis, nowadays being considered as the most famous unsolved problem in mathematics, states that the zeta function does not have non-trivial zeroes outside the line $Re\ s = \frac12 $.

	Nowadays, the independence of all events related to the divisibility of numbers, previously just discussed, can be  okay
 proved. Successful examples of using this approach enabling the same answer as a veracious probability can be found, for example, in Tenenbaum (1995), Chaumont and Yor(2003), exercise 2.10 and, especially illustrious, in Pinsky (2014). See also a fabulous book by Kac (1959) devoted to proximal topics.

Finally, the third approach to the Problem is based on the ideas of algebraic number theory. In this approach, the formulation using the {\it ring of finite integral adeles} $\widehat{\mathbb{Z}},$ which is the well-known compactification of $\mathbb{Z},$ is crucial. $\widehat{\mathbb{Z}}$ is the compact ring densely containing $\mathbb{Z}$ on which exists a unique Haar {\it probability measure} so that all tools provided by probability theory are applicable. The limit theorem solving the Problem becomes then a form of the Law of Large Numbers in the extended probability space, see Kubota and Sugita (2002), Sugita and Takanobu (2003). This approach was set up by Novoselov (1964) and developed further by many followers, see Sugita and Takanobu (2003) for their list. In the last paper, the reader can find also some promising new advancements.

\section{ A Russian perspective}

	In Russia, the Problem is often referred to as the "Chebyshev's problem" named after the outstanding Russian mathematician of the 19th century. Chebyshev, as his students reported, used to introduce it at his lectures as follows:\newline
\indent	{\it Find the probability of irreducibility of a rational fraction, the numerator and the denominator of which are chosen at random.}

	Lyapunov, probably the most famous student of Chebyshev, attended his lectures in 1879-1880 and was known to carefully record them. More than half a century later, Krylov, a notable applied mathematician and the world authority in shipbuilding mechanics, had published Chebyshev's lectures recorded by Lyapunov, see Chebyshev (1936). Translated into English, these lectures can be found in Sheynin (1991). The lectures include the problem about irreducibility of a randomly selected fraction that Chebyshev solved at the heuristic level of rigor as an example of probabilistic reasoning, see Chebyshev(1936), pp. 152-154. Because, as it follows from Lyapunov's records, Chebyshev, in his lectures, did not mention anybody's name in connection with the Problem, it seems quite natural that his students, Lyapunov included, could have referred to it as \ "Chebyshev's problem". And because Chebyshev is commonly considered as the grandfather of Russian mathematics, this term appears to be a robust mathematical folklore in Russia.
	
	Known to the authors later references to "Chebyshev's problem" in publications available in Russian include, first of all, the book by Yaglom and Yaglom (1954), see problem 90, p. 44. Note that in its English editions, Yaglom and Yaglom (1964/1987), where the Problem is included under the number 92, a reference to Chebyshev is absent. Also, the term "Chebyshev's problem" can be found in several instructional materials on Probability and Statistics, including a problem-solving book by Emelyanov and Skitovich (1967), a widely used text Sveshnikov et al.(1970), problem 4.10 with a reference the 1954 edition of Yagloms, and a book by Zhukov (2004), pp. 157-158, geared toward a popular audience. Sometimes reference to "Chebyshev problem" can be found in English sources (e.g., Erd\"{o}s and Lorentz (1959), Cacoullos (2012), problem 155. )

	Alternative evidence, along with reiteration of the Chebyshev's solution, can be found in a textbook by Markov (1913), pp. 173-176, another eminent student of the grandfather of Russian mathematics, widely known for the invention of random processes called Markov chains. Markov notes that the Problem is defined only "after a series of conditions that explain the meaning of the words that the numerator and the denominator of a fraction are chosen at random". Nonetheless, Markov almost neglects the discussion of these conditions and their role in his own reasoning.

	It should also be mentioned that heuristic solution by Chebyshev and Markov was seriously criticized by Bernstein (1964/1928) who pointed out that in the context of the Problem "one has to determine the limiting or asymptotic behavior of frequencies of natural numbers from a certain class that are distributed according to a certain rule within the sequences under consideration, rather than probabilities which we would never identify experimentally. Those limiting behaviors of ratios represent familiar analogies with mathematical probabilities, and while they are quite important for number theory from a heuristic perspective, mixing these two concepts appears to be an unfortunate mishmash" (p. 220).

\section{ Evidence leads to Dirichlet}

	In his textbook, in a footnote on p.173, Markov (1913) mentions that the Problem and its solution can also be found in the lectures by Kronecker (1901), more specifically, in Lecture 24 in which Dirichlet is given precedence for this solution over all other authors. The last remark must be elaborated further.  It is known that beginning from 1841, Kronecker studied at the University of Berlin where he attended lectures of Dirichlet and Steiner. In 1845, Kronecker defended a doctoral dissertation on algebraic number theory under the direction of Dirichlet. Soon after the defense, he left Berlin due to family circumstances and returned only in 1855, the year of Gauss's death, when Dirichlet moved to G\"{o}ttingen to take over Gauss's post. At first, Kronecker did not have an official university position, but, after being elected to Berlin Academy in 1862, started enjoying offers of lectureship appointments.

	One can suggest that in his own lectures, Kronecker (1901) reminisced Dirichlet's lectures of the 1841-1843. This brief time span preceded the period when Chebyshev lectured in St Petersburg, something that happened after he succeeded Bunyakovsky at the University of St Petersburg in 1860, see Prudnikov (1976). Thus, comparing Chebyshev and Dirichlet in the context of the Problem, one can conclude that Dirichlet encountered it earlier than Chebyshev. In addition, in 1852, from June till November, Chebyshev visited Europe and met there several distinguished mathematicians, Dirichlet included as stated in Seneta (2001). Who knows, it might be that Chebyshev learned from Dirichlet the problem about the irreducibility of a random fraction.

	In Lecture 24, Kronecker attends to the solution of the Problem. He meticulously and with great thorough carries out necessary estimates and proves that the density of pairs of relatively prime numbers over all possible pairs tends to the number $6/\pi^2.$ Then, he recalls that Dirichlet in his lectures more often considered the Problem from the probabilistic perspective. This allowed Dirichlet to reach the answer more rapidly, yet his reasoning may not be considered a rigorous proof.

	With this in mind, Kronecker analyses Dirichlet's reasoning. Let $w$ be the probability that two randomly chosen natural numbers $i$ and $k$ are relatively prime. Closely connected to $w$ is the probability $w_t$ that the greatest common factor of $i$ and $k$ is equal to $t \ge 1.$ It is clear that the original set of pairs is greater than the latter set by the factor $t^2.$ Because any pair of natural numbers has a common factor, we have $\sum_{t} w_t =1,$ whence
 $$
      w  = 1/ \sum_t t^{-2} = \frac{6}{\pi^2}.
 $$

	Further, Kronecker notes that such an approach to the problem assumes a priori that the probability in question does exist;\ that is, there exists the limiting value of the studied ratio of the number of pairs as their total number grows large and that this value is represented in an analytic form. In rigorous terms, Kronecker writes, "what Dirichlet really proves is that the probability in question, \ if after all it exists, has to be equal to \ $6/\pi^2$ ", adding that his own proof is free from this deficiency. By the way, it was Kronecker who served as editor of the posthumous publication of works by Dirichlet (1897).

	In his famous paper Dirichlet (1849), that appears also in Dirichlet (1897), Dirichlet solved a complex problem about asymptotic behavior of the sum of the values of the Euler phi function (see commentaries in Hardy and Wright (2008), p. 359) and in Dickson (1952), p. 119). Probably, this allowed Dirichlet to make his reasoning rigorous and to obtain an explicit answer for the Problem. It is quite possible though that in his lectures, Dirichlet simply did not pay attention to details. Nonetheless, Dirichlet did have keen interest in the theory of probability, a fact that a paper by Fischer (1994) supports. In particular, in this paper, a number of aspects of a lecture course on probability theory by Dirichlet have been compared to that of Chebyshev.

\section{ The problem's appeal has endured through time}

	It is worth noting, that 40 years later, Ces\`{a}ro (1881, 1883, 1884) replicated the Problem (and its solution) from Dirichlet's lectures. In particular, in the first publication, Ces\`{a}ro very briefly, literally in two lines, formulates a question about the chances of co-primality of just any two natural numbers and states that such chances are 61 to 39. Here, 0.61 is undoubtedly an approximate value of the constant $6/\pi^2$. In the subsequent publications, Ces\`{a}ro(1883, 1884), this statement was justified using arguments similar to those of Dirichlet and, in addition, some further generalizations were provided. This solution was repeated a century later in Abrams and Paris (1992), where the authors referred to Yaglom and Yaglom (1987) justifying the existence of the limit.

	In the modern textbook Bundschuh (2008), the author states that solution to the Problem was found by Ces\`{a}ro and Sylvester and then adds, "it appears that already in 1849 the solution was found by Dirichlet using somewhat different technique" (p. 52). The above-mentioned testimony of Kronecker, which, Bundschuh, most likely, did not take into account, without a doubt speaks to the favor of Dirichlet as far as the authorship of the Problem is concerned.

	The Problem did not escape interest of notable mathematicians of the 20th century, yet, for some reason, without any attention to its rich history. Indicative is a story described by Bellman (1984) in his autobiography. Shortly after the World War II, Bellman, jointly with a number theorist Shapiro, tackled the Problem and the two submitted their solution to the journal Transactions of American Mathematical Society which was edited by Kaplansky. The latter sent it to Erd\"{o}s for review who recommended acceptance. Nonetheless, Kaplansky decided to ask Erd\"{o}s if he could simplify the proof. Consequently, Erd\"{o}s found a brief and elegant proof that caused Kaplansky to suggest that the former be a co-author of the submitted paper. But Bellman and Shapiro argued that the length of their proof was not due to the deficiency of the method used but, rather, to its thoroughness and, had they skipped some details, they could have also presented a brief proof. As a result, the paper was withdrawn.

From the historical perspective, the essence of this tale is that four eminent mathematicians were not aware of the fact that the problem they discussed has a long history and that a century before them it was posed and solved by Dirichlet, well-known to Chebyshev, and, perhaps, even Gauss had been cognizant of the solution and could have shared it with Dirichlet during one of their get-togethers.

\section{ Great Gauss enters the stage}

	One can regard as an improbable speculation the above suggestion that Dirichlet himself might have learned the Problem and its solution from Gauss who, thereby, could have been the first to pose the Problem. Nonetheless, knowing that Dirichlet did meet Gauss several times (in 1827, 1828, and 1849) and enjoyed receiving his letters, see Elstrodt(2007), makes our speculation sound fairly plausible. We will return to the discussion of this hypothesis at the end of the paper.

	Furthermore, it is well known that Gauss used not to publish everything he had know\-ledge of. The most famous instances of that include non-Euclidean geometry, quaternions, and the method of least squares. There is a legend, see Choi (2011), that young and talented Jacobi (a close friend of Dirichlet) visited Gauss to show him the newest theory of elliptic functions, something that later has given him the fame. Gauss listened carefully, praised Jacobi, then opened a drawer of his table and pulled out a bunch of notes. Gauss explained to Jacobi that he completed work on this subject some time ago, but did not find the topic worthy of publishing. In turn, puzzled Jacobi asked Gauss as to why he had published even weaker results. An answer to Jacobi's question can be found at the seal of Gauss (Fig.1) where it is written in Latin "Pauca, sed matura", that is, "Few, but ripe."
\begin{figure}[h!]
\begin{center}
\includegraphics[scale=0.55]{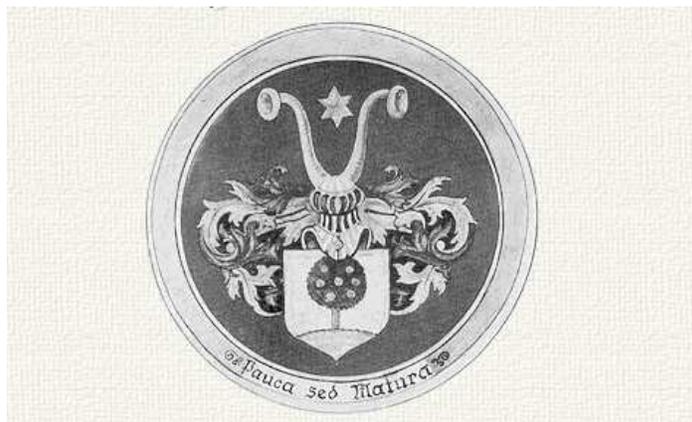}\caption{Gauss's seal}
\end{center}
\end{figure}			

\section{ Equivalent formulations and various generalizations of the Problem}

	There exist equivalent and similar formulations of the Problem worthy to be mentioned. For example, Kranakis and Pocchiola (1994), with a reference to Knuth (1981), p. 324, crediting the Problem to Dirichlet, proposed to calculate the probability that in the lattice plane, a segment connecting a randomly chosen lattice point to the origin is free from any other lattice point. It is not difficult to understand that if the segment does pass through such a point, then the fraction formed by the coordinates of the randomly chosen point is reducible. In a slightly different form, this interpretation is discussed in Apostol (1976), p. 63, Theorem 3.9 and Knuth (1981), p.324. Pieprzyk, Hardjono and Seberry (2013), pp. 190-191, narrowed down the Problem to odd integers and found that the probability of two such numbers, when randomly selected, being co-prime is equal to $8/\pi^2.$

	A natural generalization of the Problem is to find the probability/density of $k > 2$ randomly selected natural numbers being relatively prime. The answer, $1/\zeta(k),$ has been known for a long time, see Lehmer (1900) and, for a more modern presentation, Nymann (1972). A similar problem about the probability of a pairwise co-primality of $k$ randomly selected natural numbers has been solved only in the case $k = 3$ with the answer (Finch (2003), p. 110; Schroeder (2009), p. 55)
$$
Q = \frac{36}{\pi^4} \prod_{p} \left(1 - \frac{1}{(p+1)^2}\right) = 0.286747...
$$

	In the remarkable monograph on the distribution of prime numbers, an outstanding number-theorist Landau, see Landau (1909), p.69, posed the following question: what is the probability of randomly selecting a prime number? In accord with the conceptualization of the time, this question was interpreted as computing the corresponding density; that is, the limit $\lim_{x \to \infty} \pi(x)/x,$ where $\pi(x),$ as usual, denotes the
number of primes not greater than $x.$ By that time, Hadamard and Vall\'{e}e-Poussin had already proved that
$\lim_{x \to \infty} \pi(x)\ln x / x,$   exists and is equal to one, see, e.g., Hardy and Wright (1975). Thus, Landau stated that the probability sought is equal to zero. This problem, under the number 94, was included in Yaglom and Yaglom (1964/1987). Chebyshev could have also received the same result using his well-known bilateral approximation of $\pi(x)$; see, for example, Tenenbaum (1995), p. 10.

	It is of interest to compute the density/probability of randomly selecting a natural number {\it free from squares}; that is, not being divisible by any square of a natural number except one. This probability is equal to $6/\pi^2.$  It was computed by Gegenbauer (1885); see also Hardy and Wright (2008), Theorem 333.

	One may wonder, as to why the probabilities of a natural number being {\it free from squares} (squarefree) and of two natural numbers being relatively prime coincide. The answer at the physical level of rigor can be found in Schroeder (2009), p. 53. If a natural number $n$ is squarefree, then it may not be divisible by a prime number $p_k$  more than once. In other words, either $n$ is not divisible by $p_k$ , or, if is, it would not be divisible for the second time. That is,
$$
Pr (p_k^2 \nmid n ) = \left(1 - 1/p_k\right) + 1/p_k \left(1 - 1/p_k\right) = \left(1 - 1/p_k^2\right).
$$

Taking the product of these probabilities over all prime numbers (under the assumption of the independence of the corresponding events) and using identity (1) yields
$$
1/ \zeta(2) = 6/\pi^2.
$$

	Even more far-reaching generalization of the squarefree problem is to find the density of numbers free from cubes, fourth powers, and so on, up to the so-called n-free numbers; that is, numbers not divisible by any n-th power of a natural number. It appears that Gegenbauer (1885) was the first to find that this density is equal to $1/\zeta(n).$  A somewhat more modern presentation of this problem can be found in Evelyn and Linfoot(1931).

	In the paper by Hafner, Sarnak and McCurley (1993), the "probability" of relative primality of determinants of two randomly chosen  matrices with integer elements is considered. This probability, $\Delta(n)$, is computed and it is expressed in terms of the infinite product
$$
\Delta(n) = \prod_{p} \left[ 1 - \left( 1- \prod_{k=1}^n (1- p^{-k} )  \right)^2      \right].
$$

The value $\lim_{n \to \infty} \Delta(n)$  has been approximately computed by Vardi (1991), p. 174, and is equal to 0.353236... , see also Flajolet and Vardi (1996).

	It is possible to extend the Problem to other algebraic structures. Here is an example (see (Collins and Johnson, 1989) for details). The ring of Gaussian integers consists of the elements $a + bi,$ where $a$ and $b$ are integers and $i^2 =-1.$ . Suppose that two Gaussian integers are chosen at random. The probability/density that they are co-prime, in the sense taken by Collins and Johnson (1989), is equal to $6/(\pi^2 G) = 0.663700...$ , where
$$
G =\sum_{k=0}^{\infty} \frac{(-1)^k}{(2k+1)^2} = 0.915965...
$$
is known as Catalan's constant.

 In general number fields one can consider the co-primality of ideals in the ring of integer elements; the answer is expressed in terms of Dedekind zeta-function of the field, see Collins and Johnson (1989).

  An interesting generalization is situated at the confluence of analysis and number theory. Let $f(n)$ be a function "not preserving the analytical properties of $n$". The hypothesis by Erd\"{o}s and Lorentz (1959) is that for such functions
$$
  d ( \{n \in {\mathbb N} : GCD (n, \lfloor f(n) \rfloor ) =1\} ) = \frac{6}{\pi^2},  \qquad \qquad \qquad (2)
$$
where $d$ is the  asymptotic, or natural density of corresponding set of the natural numbers.
Surprisingly, the answer is again the constant $6/\pi^2$ which we met before more than once.

Erd\"{o}s and Lorentz based their argument on the earlier result by Watson (1953) stating that (2) is true if $f(n) = n\alpha,$ where $\alpha$ is irrational. They further proved the validity of (2) for differentiable and sublinear $f$, while Delmer and Deshouillers (2002) proved it when $f$  is some non-integer power of $n$ or of $\log n.$ In the recent paper by Bergelson and Richter (2016) it is proved that the hypothesis (2) is valid for $f$ from the so-called Hardy field (very broad class of differentiable functions) under some limitations on the growth of $f$, for instance, for such $f(n)$ as
$$
\log \log n, \,  \exp(n), \, \Gamma(n),\,  \zeta(n), \, Li(n),... \, .
$$

\section{Arnold's perspectives on the genesis of the Problem}

	It is interesting to recount perspectives on the Problem by Arnold (1937-2010), one of the greatest mathematical minds of the late 20th century, whose opinion bears great weight. Arnold (2003a) is very direct when referring to the problem about irreducibility of a random fraction as a ... theorem of Gauss. In Arnold (2015), it is stated that the "probability was computed by Gauss and the result published by Dirichlet" (p. 86) and Dirichlet (1849) is given as the reference to the above-mentioned result. Moreover, in Arnold (2003b), p. 4, he hypothesizes that this and similar results could have been known already to Euler, perhaps without proof, something that is not unlikely, and was later completed by Gauss.
	
In another book, that promotes the notion of mathematics as an experimental science, Arnold (2005) doubled-down on his claim about the authorship of the Problem. To this end, he considers 80 lattice points within the disk $ \{ x^2 + y^2 \le 25 \}$ (excluding the origin), finds that 48 of them have relatively prime coordinates because when such points are connected with the origin, the connecting segment is free from other lattice points, cf. Apostol (1976), p. 63; Kranakis and Pocchiola (1994), and then concludes that "the frequency of irreducibility within the circle is equal to $48/80 = 3/5$, that is, 60 \% ", Arnold (2005),  p. 13.

Arnold argues that through experimenting with larger circles one can compute the limiting probability of irreducibility which is approximately 0.608... and then goes on to suggest that "by computing this experimentally found constant, Euler obtained its exact value, $C = 6/\pi^2.$ This experimental work led him to a great deal of mathematical discoveries - theory of zeta functions, theory of Fourier series (for rough periodic functions), and to the theory of graduated algebras and their Poincar\'{e} series" (ibid, p. 14). Furthermore, Arnold believes that "it is due to Euler's experimental investigation of the probability of irreducibility of fractions that identity (1) was developed" (ibid, p. 17).

 Hence, it is possible that the author of the Problem could be also Euler, something that is in full agreement with the title of the book by Dunham (1999), dedicated to him.

From the historical perspective, all these suggestions give any indication of being incredibly fascinating and terrific, but, unfortunately, the famous author did not provide any references in support of his statements.

\section{Conclusion}

	All things considered, the authors have come to conclude as follows. Historical evidence suggests that Chebyshev was not the first to pose and heuristically solve the Problem, although, like Ces\`{a}ro, he could have done this independently. Instead, the priority should be given to Dirichlet and, perhaps, to Gauss. Besides, the Problem was revisited multiple times after Dirichlet and Chebyshev. Finally, it would be very interesting to find the traces of the Problem in the vast mathematical heritage of Gauss and Euler, because Arnold is confident that is the case. This is something that the authors, unfortunately, were not able to do.

\section{ Acknowledgement}

The authors are grateful to Dr. Galina Sinkevich for help and valuable discussions. Research of the second author was supported by grant RFBR No. 16-01-00258 and by grant SPbGU-DFG 6.65.37.2017.

\section{References}

\quad \, Abrams, A. D. and Paris, M. T. (1992). The probability that $(a, b) = 1.$ The College Mathematics Journal, 23, 47.

\medskip

Apostol, T. M. (1976). Introduction to Analytic Number Theory. Springer, New York.

\medskip

Arnold, V. I. (2003a). New Obscurantism and Russian Education. Fazis, Moscow. In Russian. Available at http://www.mccme.ru/edu/viarn/obscur.htm.

\medskip

Arnold, V. I. (2003b). Euler Groups and Arithmetic of Arithmetic Progressions. Moscow Center for Continuous Mathematics Education, Moscow.

\medskip

Arnold, V. I. (2005). Experimental Mathematics. Fazis, Moscow. In Russian.

\medskip

Arnold, V. I.  (2015). Lectures and Problems: A Gift to Young Mathematicians. American Mathematical Society, Providence, RI.

\medskip

Bellman, R. (1984). Eye of the Hurricane. An Autobiography. World Scientific, Singapore.

\medskip

Bergelson V., Richter F. K. (2016). On the density of coprime tuples of the form $(n,\lfloor f_1 (n)\rfloor,\ldots,\lfloor f_k (n)\rfloor) $, where $ f_1, \ldots, f_k $ are functions from a Hardy field. ArXiv:1611.08044.

\medskip

Bernstein, S. N. (1964/1928). Contemporary state of the theory of probability. Collected Works, vol. 4. USSR Academy of Sciences Press, Moscow. In Russian.

\medskip

Bundschuh, P. (2008). Einfuhrung in die Zahlentheorie, 6 Auflage, Springer-Verlag. In German.

\medskip

Cacoullos, T. (2012). Exercises in probability.  Springer Science \& Business Media.

\medskip

Ces\`{a}ro, E. (1881). Question propos\'{e}e 75. Mathesis 1, 184.

\medskip

Ces\`{a}ro, E. (1883). Question 75 (Solution). Mathesis 3, 224-225.

\medskip

Ces\`{a}ro, E. (1884). Probabilit\'{e} de certains faits arithm\'{e}tiques. Mathesis 4, 150-151.

\medskip

Chaumont, L., Yor, M. (2003). Exercises in Probability, Cambridge University Press.

\medskip

Chebyshev, P. L. (1936). Theory of Probability. Lectures delivered in 1879-1880. Published from notes taken by A. M. Lyapunov (A. N. Krylov, ed.). USSR Academy of Sciences Press, Moscow-Leningrad. In Russian.

\medskip

Choi, Y. (2011, Jan. 25). Re: Jacobi and Gauss. [Blog comment]. Retrieved from https://ifwisdomwereteachable.wordpress.com/2011/01/25/lesprit-descalier/

\medskip

Collins, G. E. and Johnson, J. R. (1989). The probability of relative primality of Gaussian integers. In Symbolic and Algebraic Computation, Proceedings of International Symposium ISSAC 1988 (P. Gianni, ed.), 252-258. Lecture Notes in Computer Science 358. Springer, New York.

\medskip

Delmer F., Deshouillers J. M.(2002). On the probability that $n$ and $[n^c]$ are coprime. Period. Math. Hung., 45,  15 -- 20.

\medskip

Dickson, L. E. (1952). History of the Theory of Numbers, vol. 1: Divisibility and Primality. New York: Chelsea.

\medskip

Dirichlet, P. (1849). Uber die Bestimmung der mittleren Werte in der Zahlentheorie. Abhandl. Kgl. Preuss. Acad. Wiss., Berlin., 69-83.

\medskip

Dirichlet, P. (1897). Werke, Bd. 2. Berlin,  60-64.

\medskip

Dunham W. (1999). Euler: The master of us all. Washington D.C.: Mathematical Association of America.

\medskip

Elstrodt, J. (2007). The life and work of Gustav Lejeune Dirichlet (1805-1859). Analytic Number Theory, Clay Mathematics Proceedings 7, 1-37. Providence, RI: American Mathematical Society.

\medskip

Emelyanov, G. V. and Skitovich, V. P. (1967). Problem-Solving Book on Probability and Statistics. Leningrad University Press, Leningrad. In Russian.

\medskip

Erd\"{o}s, P., Lorentz, G.G. (1959). On the probability that $n$ and $g(n)$ are relatively prime. Acta Arith., 5, 35 -– 44.

\medskip

Evelyn, C. J. A. and Linfoot, E. H. (1931). On a problem in the additive theory of numbers. Annals of Mathematics 32, 261-270.

\medskip

Finch, S. R. (2003). Mathematical Constants. Cambridge University Press, New York.

\medskip

Fischer, H. (1994). Dirichlet's contributions to mathematical probability theory. Historia Mathematica 21, 39-63.

\medskip

Flajolet, P. and Vardi, I. (1996). Zeta Function Expansions of Classical Constants. Available at \\ http://algo.inria.fr/flajolet/Publications/FlVa96.pdf

\medskip

Gegenbauer, L. (1885). Asymptotische Gesetze der Zachlentheorie. Denkschriften Akad. Wiss. Wien, 49, 37-80.

\medskip

Hafner, J. L., Sarnak, P. and McCurley, K. (1993). Relatively prime values of polynomials. In Tribute to Emil Grosswald: Number Theory and Related Analysis (M. Knopp and M. A. Seingorn, eds.) 437-444. American Mathematical Society, Providence, RI.

\medskip

Hardy, G. H. and Wright, E. M. (1975). An Introduction to the Theory of Numbers (4th edition). Oxford University Press, Oxford.

\medskip

Hardy, G. H. and Wright, E. M. (2008). An Introduction to the Theory of Numbers (6th edition). Revised by D. R. Heath-Brown and J. H. Silverman. With a foreword by A. Wiles. Posts and Telecommunications Press, Beijing.

\medskip

Hombas, V. (2013). What's the probability of a rational ratio being irreducible? International Journal of Mathematical Education in Science and Technology 44, 408-410.

\medskip

Kac, M. (1959). Statistical Independence in Probability, Analysis and Number Theory. The Carus Mathematical Monographs, N 12. The Mathematical Association of America, Rahway, NJ.

\medskip

Knuth, D. E. (1981). The Art of Computer Programming: Seminumerical Algorithms, v. 2. Addison-Wesley, Reading, MA.

\medskip

Kranakis, E. and Pocchiola, M. (1994). Camera placement in integer lattices. Discrete and Computational Geometry 12, 91-104.

\medskip

Kronecker, L. (1901). Vorlesungen ueber Mathematik, Teil 2, Abschnitt 1 (Vorlesungen ueber Zahlentheorie), Bd. 1. Springer, Berlin Heidelberg GmbH. In German.

\medskip

Kubota, H. and Sugita, H.(2002). Probabilistic proof of limit theorems in number theory by means of adeles. Kyushu Journal of Mathematics, 56(2), 391--404.

\medskip

Landau, E. (1909). Handbuch der Lehre von der Verteilung der Primzahlen. BG Teubner, Leipzig und Berlin. In German.

\medskip

Lehmer, D. N. (1900). Asymptotic evaluation of certain totient sums. American Journal of Mathematics 22, 293-335.

\medskip

Markov, A. A. (1913). Calculus of Probabilities (3rd edition). Imperial Academy of Sciences Press, St Petersburg. In Russian.

\medskip

Mertens, F. (1874). Ueber einige asymptotische Gesetze der Zahlentheorie. J. reine und angew. Mathem. 77, 289-338.

\medskip

Novoselov E. V. (1964). A new method in probabilistic number theory. Izv. Ross. Akad. Nauk. Ser. Matem. 28(2), 307--364 (in Russian).

\medskip

Nymann, J. E. (1972). On the probability that k positive integers are relatively prime. Journal of Number Theory 4, 469-473.

\medskip

Pieprzyk, J., Hardjono, T. and Seberry, J. (2013). Fundamentals of Computer Security. Springer, Berlin.

\medskip

Pinsky, R. G. (2014). Problems from the Discrete to the Continuous. Probability, Number Theory, Graph Theory, and Combinatorics. Springer, New York.

\medskip

Prudnikov, V. E. (1976). Pafnuty Lvovich Chebyshev, 1821-1894. Nauka, Moscow. In Russian.

\medskip

Schroeder, M. (2009). Number Theory in Science and Communication: With Applications in Cryptography, Physics, Digital Information, Computing, and Self-Similarity. Springer, New York.

\medskip

Seneta, E. (2001). Pafnutii Lvovich Chebyshev (or Tchebichef). Statisticians of the Centuries, Heyde C.C., Seneta E., Crepel, P., Fienberg, S.E., Gani, J. (Eds.),  Springer, New York, 176-180.

\medskip

Sheynin, O. (1994). Chebyshev's lectures on the theory of probability. Archive for History of Exact Sciences, 46(4), 321-340.

\medskip

Sugita, H. and Takanobu, S. (2003). The probability of two integers to be co-prime, revisited —- on the behavior of CLT-scaling limit. Osaka J. Math, 40(4), 945-976.

\medskip

Sveshnikov, A. A., Ganin, M. P., Diner, I. Y, Komarov, L. B., Starobin, K. B. and Volodin, B. G. (1970). Problem-Solving Book on Probability, Statistics, and Random Functions. Nauka, Moscow. In Russian.

\medskip

Sylvester, J. J. (1883). Sur le nombre de fractions ordinaires in\'{e}gales qu'on peut exprimer en se servant de chiffres qui n'exc\`{e}dent pas un nombre donn\'{e}. C. R. Acad. Sci. Paris XCVI, 409-413. Reprinted in H. F. Baker (Ed.), The Collected Mathematical Papers of James Joseph Sylvester, vol. 4, Cambridge University Press, 86.

\medskip

Tenenbaum, G. (1995). Introduction to analytic and probabilistic number theory. Cambridge University Press, New York.

\medskip

Vardi, I. (1991). Computational Recreations in Mathematica. Addison-Wesley, Boston.

\medskip

Watson G. L.(1953)  On integers $n$ relatively prime to $[\alpha n ]$. Canadian J. Math. 5,  451 -– 455.

\medskip

Yaglom, A. M. and Yaglom, I. M. (1954/1964/1987). Challenging Mathematical Problems with Elementary Solutions, vol. 1: Combinatorial Analysis and Probability Theory. Dover, New York. (Russian edition 1954).

\medskip

Zhukov, A. V. (2004). Ubiquitous Number $ \pi $. Editorial URSS, Moscow. In Russian.

\end{document}